\documentstyle{amsppt}
\tolerance 3000
\pagewidth{5.5in}
\vsize7.0in
\magnification=\magstep1
\NoBlackBoxes
\NoRunningHeads
\widestnumber \key{AAAAAAA}
\topmatter
\author A. Iosevich and I. {\L}aba
\endauthor
\title Discrete subsets of ${\Bbb R}^2$ and the associated distance 
sets
\endtitle
\date March 14, 2002
\enddate
\abstract We prove that a well-distributed subset of
${\Bbb R}^2$ can have a separated distance set only if the
distance is induced by a polygon.
\endabstract
\thanks Research of A. Iosevich supported in part by the NSF
Grant DMS00-87339 \endthanks
\thanks Research of I. {\L}aba supported in part by the NSERC
Grant 22R80520
\endthanks
\address A. Iosevich, University of Missouri
\ email: iosevich \@ wolff.math.missouri.edu
\endaddress
\address I. {\L}aba, University of British Columbia
\ email: ilaba \@math.ubc.ca
\endaddress
\endtopmatter

\head Basic definitions \endhead

\vskip.125in

\definition{Separated sets} We say that $S \subset {\Bbb R}^d$ is 
separated
if there exists $c_{separation}>0$ such that $||x-y|| \ge
c_{separation}$ for every $x,y\in S$, $x\neq y$. Here and throughout 
the paper,
$||x||=\sqrt{x_1^2+\dots+x_d^2}$ is the standard Euclidean distance.
\enddefinition

\definition{Well-distributed sets} We say that $S \subset {\Bbb R}^d$
is well-distributed if there exists a $C_{density}>0$ such that every
cube of side-length $C_{density}$ contains at least one element of $S$.
\enddefinition

\definition{$K$-distance} Let $K$ be a bounded convex set, symmetric
with respect to the origin. Given $x,y \in {\Bbb R}^d$, define
the $K$-distance, ${||x-y||}_K=\inf \{t: x-y \in tK\}$.
\enddefinition

\definition{$K$-distance sets} Let $A \subset {\Bbb R}^2$. 
Define $\Delta_K(A)=\{{||x-y||}_K: x,y \in A\}$, the $K$-distance set 
of $A$.
\enddefinition

Notation: $A \lessapprox B$, with respect to the parameter $R$,
means that there exists a positive constant $C_{\epsilon}$ such
that $A \leq C_\epsilon R^{\epsilon}B$ for any $\epsilon>0$. Similarly, 
$A
\lesssim B$ means that there exists a $C>0$ such that $A \leq CB$,
and $A \approx B$ means that $A \lesssim B$ and $B \lesssim A$.

\head Introduction \endhead

\vskip.125in

Distance sets play an important role in combinatorics and its
applications to analysis and other areas. See, for example,
\cite{AP95} and the references contained therein. Perhaps the
most celebrated classical example is the Erd\H{o}s Distance
Problem, which asks for the smallest possible cardinality of
$\Delta_{B_2}(A)$, where $A\subset{\Bbb R}^2$ has cardinality
$N<\infty$ and $B_2$ is the Euclidean
unit disc. Erd\H{o}s conjectured that $\# \Delta_{B_2}(A) 
\gtrapprox N$. The best known result to date in
two dimensions is due to Solymosi and T\'oth who prove in
\cite{ST01} that $\# \Delta_{B_2}(A) \gtrsim N^{\frac{6}{7}}$.
For a survey of higher dimensional results see \cite{AP95} and
the references contained therein. For applications of distance
sets in analysis, see, for example, \cite{IKP99}, where distance
sets are used to study the question of existence of orthogonal
exponential bases.

The situation changes drastically if the Euclidean disc $B_2$ is
replaced by a convex planar set with a ``flat" boundary. For example,
suppose that $K={[-1,1]}^2$, which corresponds to the ``taxi-cab"
or $l^1({\Bbb R}^2)$ distance. Let $A=\{m \in {\Bbb Z}^2: 0 \leq
m_i \leq N^{\frac{1}{2}}\}$. Then $\# A \approx N$, and it is
easy to see that $\# \Delta_K(A) \approx N^{\frac{1}{2}}$, which
is much less than what is known to be true for the Euclidean
distance. In fact, it follows from an argument due to Erd\H{o}s
(\cite{Erd46}; see also \cite{I01}) that the estimate $\# \Delta_K(A)
\gtrsim  N^{\frac{1}{2}}$ holds for any $K$.

The example in the previous paragraph shows that the properties
of the the distance set very much depend on the underlying
distance. One way of bringing this idea into sharper focus is
the following. Let $S$ be a separated subset of ${\Bbb R}^2$,
$\alpha$-dimensional in the sense that
$$ \# (S \cap {[-R,R]}^2) \approx R^{\alpha}. \tag0.1$$

If Erd\H{o}s' conjecture holds, then $\# \Delta_{B_2}(S \cap
{[-R,R]}^2) \gtrapprox R^{\frac{2 \alpha}{2}}$. This would imply
that if $\alpha>1$ then $\Delta_{B_2}(S)$ is not separated. This
formulation expresses the Erd\H{o}s Distance Conjecture in the
language of the Falconer Distance Conjecture (see e.g.
\cite{Wolff02}) which says that if a compact set $E \subset
{\Bbb R}^2$ has Hausdorff dimension $\alpha>1$, then
$\Delta_{B_2}(E)$ has positive Lebesgue measure.

On the other hand, as we have seen above, if $K={[0,1]}^2$, $S$
may be $2$-dimensional (for example, the integer lattice), and
the distance set $\Delta_K(S)$ is, nevertheless, separated. The
purpose of this paper is to address the following question.
Suppose that $S$ is a well-distributed subset of
${\Bbb R}^2$. For which $K$'s can $\Delta_K(S)$ be separated? We
conjecture that this can only happen if $K$ is a polyhedron.
This is indeed what we prove in two dimensions.

\proclaim{Theorem 0.1} Let $S$ be a well-distributed
subset of ${\Bbb R}^2$. Suppose that $\Delta_K(S)$ is separated.
Then $K$ is a polygon with finitely many sides. \endproclaim

\vskip.125in

\head Proof of Theorem 0.1 \endhead

\vskip.125in

Let $S \subset {\Bbb R}^2$ be a well-distributed set. We may
assume that $C_{density}=1$. We identify
${\Bbb R}^2$ with ${\Bbb C}$ via $z=re^{i\theta}$, and denote by
$C_{\theta_1,\theta_2}$ the cone $\{re^{i\theta}:\
\theta_1<\theta<\theta_2\}$. We also write $\Gamma=\partial K$.
A {\it line segment} will always be assumed to have non-zero
length.

Theorem 0.1 is an immediate consequence of Lemmas 1.1 and 1.2 below.

\proclaim{Lemma 1.1} Let $S$ be as above. If $\Delta_K(S)$ is
separated, then for any $\theta_1 <\theta_2$ the curve
$\Gamma\cap C_{\theta_1,\theta_2}$ contains a line segment.
\endproclaim

\proclaim{Lemma 1.2} Let $S$ be as above. If $\Delta_K(S)$ is
separated, then $\Gamma$ may contain only a finite number of
line segments such that no two of them lie on one straight line.
\endproclaim

We now prove Lemmas 1.1 and 1.2. The main geometrical
observation is contained in the next lemma.

\proclaim{Lemma 1.3} Let $\Gamma=\partial K$, where $K \subset
{\Bbb R}^2$ is convex. Let $\alpha>0$, $x \in {\Bbb R}^2$,
$x\neq 0$.

(i) If $\Gamma\cap (\alpha\Gamma+x)$ contains three distinct
points, at least one of these points must lie on a line segment
contained in $\Gamma$.

(ii) $\Gamma\cap (\alpha\Gamma+x)$ cannot contain more than 2 line
segments such that no two of them lie on one line.
\endproclaim

We will first prove Lemmas 1.1 and 1.2, assuming Lemma 1.3;
the proof uses a variation on an argument of Moser \cite{Mo}.
The proof of Lemma 1.3 will be given later in this section.

Fix 2 points $P,Q\in S$; translating $S$ if
necessary, we may assume that $P=-Q$. Let
$$A_N=\{z\in {\Bbb R}^2:\ d_K(x,0)\in (10N,10(N+1))\}. \tag1.1$$
Observe that for all $N$ large enough (depending on $\theta_1,
\theta_2$) the set $A_N\cap C_{\theta_1,\theta_2}$ contains at least
$N(\theta_2-\theta_1)$ points of $S$.

Observe that for all $s\in S\cap A_N$ we have
$$
d_K(s,P)\in[10N-d_K(P,0),10(N+1)+d_K(Q,0)],$$
and similarly for $d_K(s,Q)$.  The length of the interval
above is independent of $N$. Since $\Delta_K(S)$ is separated,
it follows that the number of distinct distances
$$ d_K(s,P),\ d_K(s,Q):\ s\in S\cap A_N \tag1.2$$
is bounded by a constant $L$ uniformly in $N$.

\subhead{Proof of Lemma 1.1} \endsubhead We may assume that $0<
\theta_2-\theta_1 < \pi/2$.  Fix $\theta'_1,\theta'_2$ so that
$\theta_1<\theta'_1<\theta'_2 <\theta_2$, and let
$C=C_{\theta_1,\theta_2}$, $C'=C_{\theta'_1,\theta'_2}$. Then
for all $N$ large enough we have
$$
C'\cap A_N \subset (C+P)\cap(C+Q) \tag1.3$$ and
$$
\#(S\cap C'\cap A_N)\geq N(\theta'_2-\theta'_1). \tag1.4$$

Let
$$
\{d_1,\dots,d_l\}=\{d_K(s,P):\ s\in S\cap A_N\},\
\{d'_1,\dots,d'_{l'}\}=\{d_K(s,Q):\ s\in S\cap A_N\}, \tag1.5$$
where $d_j\neq d_k$, $d'_j\neq d'_k$ if $j\neq k$.
Then $l,l'\leq L$ with $L$ independent of $N$ (see (1.2)). We have
$$
S\cap C'\cap A_N\subset S\cap C'\cap \bigcup_{i,j}\Gamma_i\cap 
\Gamma'_j,
\tag1.6$$
where $\Gamma_i=d_i\Gamma+P$, $\Gamma'_j=d'_j\Gamma+Q$. Thus if $N$ is
large enough ($N\geq 10L^2(\theta'_2-\theta'_1)^{-1}$ suffices),
there are $i,j$ such that
$\#(S\cap C'\cap\Gamma_i\cap\Gamma'_j)\geq 10$.  It follows from Lemma
2.3(i) that at least one of the points in $S\cap C'\cap\Gamma_i$ lies
on a line segment $I$ contained in $C'\cap\Gamma_i$. By (1.3),
$I$ is contained in $d_i\Gamma\cap C$, hence $\Gamma\cap C$ contains
the line segment $d_i^{-1}I$.

\subhead{Proof of Lemma 1.2} \endsubhead Suppose that $\Gamma$
contains line segments $I_1,\dots,I_M$, all pointing in different
directions.

We will essentially continue to use the notation of the proof of
Lemma 1.1. Choose $P,Q$ as above, and let $C_m$ denote cones
$C_m=C_{\theta_m,\theta'_m}$ such that $\theta_m<\theta'_m$ and
$\Gamma\cap C_m\subset I_m$. Let also $C'_m\subset C_m$ be
slightly smaller cones. Let $N$ be large enough so that each
sector $C'_m\cap A_N$ contains at least $10$ points of $S$, and
so that
$$ A_N\cap C'_m\subset (C_m+P)\cap(C_m+Q). \tag1.7
$$
Let also $d_i,d'_j, \Gamma_i,\Gamma'_j$ be as above.  Then
for each $m$
$$
S\cap A_N\cap C'_m\subset \bigcup_{i,j}\Gamma_i\cap
\Gamma'_j\cap C'_m. \tag1.8$$
If $N$ is large enough, $\Gamma_i\cap C'_m\subset C_m$
and $\Gamma'_j\cap C'_m\subset C_m$, hence
the set on the right is a union of line segments parallel to
$I_m$.  It must contain at least one such segment, since the
set on the left is assumed to be non-empty.  Therefore the
set
$$
\bigcup_{i,j}\Gamma_i\cap\Gamma'_j \tag1.9
$$
contains at least $M$ line segments pointing in different directions,
one for each $m$.  But on the other hand,
by Lemma 2.3(ii) any $\Gamma_i\cap \Gamma'_j$ can
contain at most two line segments that do not lie on one line.
It follows that the set in (1.9) contains at most $2L^2$ line
segments in different directions, hence $M\leq 2L^2$.  Since
$\Gamma$ can contain at most two parallel line segments that
do not lie on one line, the number of line
segments in Lemma 1.2 is bounded by $4L^2$ as claimed.

\subhead{Proof of Lemma 1.3} \endsubhead We first prove part (i)
of the lemma.  Suppose that $P_1,P_2,P_3$ are three distinct
points in $\Gamma\cap(\alpha\Gamma +x)$.  We may assume that
they are not collinear, since otherwise the conclusion of the
lemma is obvious.  We have $P_1,P_2,P_3\in\Gamma$ and
$P'_1,P'_2,P'_3\in\Gamma $, where $P'_j=\alpha^{-1}(P_j-x)$. Let
$T$ and $T'$ denote the triangles $P_1P_2P_3$ and
$P'_1P'_2P'_3$, and let $K'$ be the convex hull of $T\cup T'$.
Since $K'\subset K$ and all of the points $P_j,P'_j$ lie on
$\Gamma=\partial K$, they must also lie on $\partial K'$.

Observe that $\partial K'$ consists of some number of the edges
of the triangles $T,T'$ and at most 2 additional line segments.
If $\partial K'$ contains at most one of $P_iP_j$ and $P'_iP'_j$
for each $i,j$, $K'$ is a polygon with at most 5 edges, hence at
most 5 vertices.  Thus if the 6 points $P_j,P'_j$ lie on
$\partial K'$, at least three of them must be collinear, and
one of them must be $P_j$ for some $j$ (otherwise the $P'_j$
would be collinear).  If these three points are distinct,
then $\Gamma$ contains the line segment joining all of them,
and we are done.  Suppose therefore that they are not
distinct.  It suffices to consider the cases when
$P_1=P'_1$ or $P_1=P'_2$.  If $P_1=P'_1$, then we must have
$\alpha\neq 1$ and $P_1,P_2,P'_2$ are distinct and collinear; if
$P_1=P'_2$, then $P'_1,P_1,P_2$ are distinct and collinear.
Thus at least three of the points $P_j,P'_j$ are distinct and
collinear, and we argue as above.

It remains to consider the case when
$\partial K'$ contains both $P_iP_j$ and
$P'_iP'_j$ for some $i,j$.  The outward unit normal vector to
$P_iP_j$ and $P'_iP'_j$ is the same, hence all four points
$P_i,P_j,P'_i,P'_j$ are collinear, at least three of them are
distinct, and $\Gamma$ contains a line segment joining all of
them.

Part (ii) of the lemma is an immediate consequence of the
following. Let $(x_1,x_2)$ denote the rectangular coordinates in
the plane.

\proclaim{Lemma 1.4} Let $I$ be a line segment contained in
$\Gamma\cap(\alpha\Gamma+u)$, where $u=(c,0)$.

(i) If $\alpha=1$, then $I$ is parallel to the $x_1$-axis.

(ii) If $\alpha\neq 1$, then the straight line containing $I$
goes through the point $\left(\frac{c}{1-\alpha},0 \right)$.
\endproclaim

\demo{Proof of Lemma 1.4} Part (i) is obvious; we prove (ii). If
$I$ lies on the line $x_2=ax_1+b$, then so does $\alpha I + u$.
But on the other hand $\alpha I + u$ lies on the line

$$ x_2=\alpha \left(a\frac{x_1-c}{\alpha}+b \right)=
ax_1-ac+\alpha b. \tag1.10$$

It follows that $b=\alpha b-ca$, hence $-\frac{b}{a}
=\frac{c}{1-\alpha}$. But $-{b}/{a}$ is the $x_1$-intercept of
the line in question.

Similarly, if $I$ lies on the line $x_1=b$, then $\alpha I + u$
lies on the lines $x_1=b$ and $x_1=\alpha b+c$, hence
$b=\frac{c}{1-\alpha}$. \enddemo

To finish the proof of Lemma 2.3 (ii), it suffices to observe
that in both of the cases (i), (ii) of Lemma 2.4 the boundary of
a convex body cannot contain three such line segments if no two
of them lie on one line.

\newpage

\head References \endhead

\vskip.125in

\ref \key AgPa95 \by P. Agarwal and J. Pach \paper Combinatorial
Geometry \yr 1995 \jour Wiley-Interscience Series \endref

\ref \key Erdos46 \by P. Erd\H{o}s \paper On sets of distances of
$n$ points \yr 1946 \jour Amer. Math. Monthly \vol 53 \pages
248-250 \endref

\ref \key I01 \by A. Iosevich \paper Curvature, combinatorics
and the Fourier transform \yr 2001 \jour Notices of the AMS \vol
46 no. 6 \pages 577-583 \endref

\ref \key IKP99 \by A. Iosevich, N. Katz, and S. Pedersen \paper
Fourier basis and the Erd\H os distance problem \jour Math.
Research Letter \yr 1999 \vol 6 \pages \endref

\ref \key Mo \by L. Moser \paper On the different distances 
determined by $n$ points \yr 1952 \jour Amer. Math. Monthly
\vol 59 \pages 85-91 \endref

\ref \key ST01 \by J. Solymosi and Cs. D. T\'oth \paper Distinct
distances in the plane \jour Discrete Comput. Geom. \vol 25 \yr
2001 \pages 629-634 \endref

\ref \key Wolff02 \by T. Wolff \paper Lecture notes in harmonic
analysis (revised) \yr 2002 \endref

\vskip.25in

\enddocument